\documentclass[12pt]{amsart}
\usepackage{amsmath,amssymb,amsfonts,amsthm}
\usepackage{mathrsfs}
\usepackage[all]{xy}
\usepackage{stmaryrd}
\usepackage{mathtools}
\usepackage{indentfirst}
\usepackage{verbatim}
\usepackage{hyperref}
\usepackage[normalem]{ulem}
\usepackage{tikz-cd}
\usepackage[shortlabels]{enumitem}
\usepackage{comment}

\newtheorem*{thm*}{Theorem}
\newtheorem{thm}{Theorem}[section]
\newtheorem{prop}[thm]{Proposition}
\newtheorem{lemma}[thm]{Lemma}

\theoremstyle{definition}

\newtheorem{defn}[thm]{Definition}
\newtheorem{rmk}[thm]{Remark}
\newtheorem{ex}[thm]{Example}

\newcommand{\on}{\operatorname}

\usepackage[margin=0.75in]{geometry}

\newcommand\cE{{\mathcal{E}}}
\newcommand\cF{{\mathcal{F}}}
\newcommand\cG{{\mathcal{G}}}

\newcommand\cI{{\mathcal{I}}}

\newcommand\cO{{\mathcal{O}}}

\newcommand\bC{{\mathbb C}}

\newcommand\bF{{\mathbb F}}

\newcommand\bN{{\mathbb N}}
\newcommand\bP{{\mathbb P}}
\newcommand\bQ{{\mathbb Q}}
\newcommand\bR{{\mathbb R}}

\title{Existence of complements for foliations}
\author{Yen-An Chen, Dongchen Jiao, and Pascale Voegtli}
\address{Department of Mathematics, Imperial College London, 180 Queen’s Gate, London SW7 2AZ, UK}
\email{yen-an.chen@imperial.ac.uk}
\address{Department of Mathematics, Brunel University London UB8 3PH}
\email{dongchen.jiao@brunel.ac.uk}
\address{Department of Mathematics, University College London, WC1E 6BT}
\email{pascale.voegtli.20@ucl.ac.uk}

\begin{document}
\maketitle
\begin{abstract}
    This paper demonstrates the existence of $\mathbb{Q}$-complements for algebraically integrable log-Fano foliations on klt ambient varieties. Additionally, we investigate properties of algebraically integrable Fano foliations such as a partial inversion of  adjunction as well as a connectedness principle. 
\end{abstract}
\section{Introduction}
A monotone  $\mathbb{Q}$-complement of a log canonical pair $(X,B)$ is a $\mathbb{Q}$-divisor $B^{+} \geq B$ such that ${K_X + B^{+} \sim_{\mathbb{Q}} 0}$ and $(X,B^{+})$ is log canonical.
The investigation of complements of Fano varieties has a rather long history. Studied in full generality for the first time by Prokhorov and Shokurov in \cite{Prokhorov_2001} the concept dates back to a publication of Shokurov \cite{shokurov2001complements} on low dimensional varieties. In 2016 in \cite{Antipluri} Birkar proved the existence of monotone $n$-complements for Fano type varieties of given dimension under the additional assumption that the coefficients of the $\mathbb{Q}$-divisor $B$ lie in a controlled set. More precisely, he shows that for all log canonical  Fano type pairs $(X,B)$ where the dimension of $X$ as well as the coefficient set of $B$ are fixed, there exists an $n$ only depending on the latter two conditions such that  ${n(K_X + B^{+}) \sim 0}$. These are commonly referred to as bounded complements. The existence of bounded complements allowed Birkar in a subsequent publication \cite{SingofLinsyst} to verify a long standing conjecture due to Borisov-Alexeev-Borisov (BAB) which states that the set of $\epsilon$-lc Fano varieties of a fixed dimension forms a bounded family.

Since the release of \cite{Antipluri}, complements have become a rather active field of research as a plethora of publications on complements testifies. Without any claim of completeness, this includes complements in particular settings such as for Fano fibrations \cite{choi2024existence} or klt-complements over lower dimensional varieties such as \cite{chen2024boundedness} or in the realm of generalized pairs such as \cite{chen2020boundedness1} and \cite{chen2020boundedness} or finally more explicit results such as \cite{figueroa2022complements} finding effective bounds on the index of $n$-complements of given coregularity.

Recent advances in the development of the Minimal Model Program (MMP) for foliated varieties, in particular on algebraically integrable foliations (see \cite{cascini2023mmp} and \cite{LMX2024minimal}) allow to investigate the existence of (bounded) complements in the foliated setting. As opposed to complements for varieties, even the existence of $\mathbb{Q}$-complements due to the general failure of Bertini-type results is by no means trivial and the main content of the present publication. By a $\mathbb{Q}$-complement for a foliated triple $(X,\mathcal{F},B)$ we mean a $\mathbb{Q}$-divisor $B^{+} \geq B$ such that ${K_{\mathcal{F}} + B^{+} \sim_{\mathbb{Q}} 0}$ and $(\cF,B^{+})$ is log canonical. So far we are unfortunately unable to demonstrate the existence of bounded complements in full generality. Partially because we fall back on classical vanishing theorems to compensate for the lack of their analogues for $K_{\cF}$ in the process of lifting sections from invariant subvarieties to $X$, partially because we face technical difficulties when trying to settle the base case for on the closure of general fibres. However, we are confident to overcome these difficulties in ongoing work, in which we moreover address boundedness questions for Fano foliations. We see our publication as a partial resolution of a problem set out by Cascini and Spicer at a conference \cite{conjecture}. Our main theorem on the existence of $\mathbb{Q}$-complements is the following:
\begin{thm}[{$=$ Theorem \ref{Qcom}}]
    Let $(X,\mathcal{F},B)$ be an lc Fano foliated triple with $\cF$ algebraically integrable. 
    Then there exists a $\bQ$-complement. 
    That is, there exists an effective divisor $D\sim_\bQ-(K_{\mathcal{F}}+B)$ such that $(X,\mathcal{F},B+D)$ is log canonical.
\end{thm}
As mentioned in the abstract, we furthermore provide proofs of a number of side results on Fano and log Calabi Yau foliations that may be of independent interest. To be quoted first, a weak form of a foliated inversion of a adjunction to foliation invariant divisors:
\begin{thm}[{$=$ Theorem \ref{inv}, Inversion of adjunction}]
    Let $(X,\mathcal{F},B)$ be a foliated triple with $\mathcal{F}$ being algebraically integrable. 
    Let $S\subseteq X$ be a prime invariant divisor with normalisation $\nu\colon S^\nu\to S$ such that $(S^\nu,\mathcal{F}_{S^\nu},\on{Diff}_{S^\nu}(\mathcal{F},B))$ is log canonical. Then
    \[\on{Nlc}(X,\mathcal{F},B)\cap S=\varnothing.\]
\end{thm}
It is to be pointed out that we exclusively deduce log canonicity of a pair from log canonicity of its restriction to an invariant divisor, however cannot prove the preservation of other types of singularities.
In the same spirit, we have:
\begin{prop}[{$=$ Proposition \ref{horizofNlc}}]
Let $(Y,\mathcal{G},\tilde{B})$ be a foliated triple where $\mathcal{G}$ is induced by a dominant rational map to a normal variety $Z_Y$ and such that $K_{\mathcal{G}}+\tilde{B} \sim_{\mathbb{Q}}0$. Then $\on{Nlc}(Y,\mathcal{G},\tilde{B})$ dominates $Z_Y$.
\end{prop}

This can be seen as a weak inversion of subadjunction for foliated log Calabi-Yau triples. The ``weak'' is to be understood in the same way as above. We further notice that even though we phrase the statement here only for restrictions to fibres, it equally applies to general members of covering families of invariant subvarieties of any dimension. 

Lastly, we include the following foliated version of a connectedness principle for algebraically integrable log Fano foliations:
\begin{lemma}[{$=$ Lemma \ref{conlem}}]
    Let $(\mathcal{F},\Delta)$ be a foliated pair on a normal variety $X$ such that
    \begin{enumerate}
        \item $-(K_{\mathcal{F}}+\Delta)$ is big and nef;
        \item $(\mathcal{F},\Delta)$ is log canonical;
        \item $\mathcal{F}$ is induced by a morphism $g\colon X\to Z$.
    \end{enumerate}
    Then $\lfloor\Delta\rfloor$ is connected. 
\end{lemma}

\section*{Acknowledgements}
The authors would like to thank Imperial College London for the wonderful research environment. 
They would also like to express their gratitude to Paolo Cascini, Anne-Sophie Kaloghiros, Ching-Jui Lai, and Calum Spicer for helpful discussions.
The second author is supported by EPSRC DTP studentship (EP/V520196/1) Brunel University London and would like to show appreciation to his supervisor Anne-Sophie Kaloghiros for useful discussion and invaluable suggestions.\par
The third author would like to thank the LSGNT for laying the mathematical foundations of this work by providing a stimulating and supporting working environment.
He was supported by the Engineering and Physical Sciences Research Council [EP/S021590/1]. 
The EPSRC Centre for Doctoral Training in Geometry and Number Theory (The London School of Geometry and Number Theory), University College London.

\section{Sketch of the proof strategy}
Analogous to the classical case we aim at an inductive proof strategy, this is we seek to invoke the existence of foliated complements assuming their existence on invariant subvarieties.\par
First, following the recent work on foliated adjunction by Cascini and Spicer \cite{cascini2023foliation}, we prove lc adjunction and inversion of adjunction for algebraically integrable foliations, which allows us to apply adjunction to invariant subvarieties to get a Fano triples in lower dimensions.\par
Second, with the recently developed tool on foliated minimal model program, we prove that the generic log canonical property defined in \cite[Definition 3.8]{araujo2014fano} coincides with McQuillan's definition \cite[Definition I.1.5]{mcquillan2008canonical} for global log canonical under suitable condition.\par
Last, since Bertini's Theorem \cite[Lemma 4.7.1]{1996alg.geom..1026K} provides us with $\mathbb{Q}$-complements naturally for Fano pairs, we may lift them to global sections and argue that they are log canonical with the above theorem in the second step.\par
Moreover, we would like to point out that since we do not have a foliated version of Kawamata-Viehweg vanishing to replace Serre vanishing, this method does not guarantee the existence of $n$-complements for bounded $n$.  

\section{Preliminaries}
We work over an algebraically closed field of characteristic zero.
\begin{defn}
    Let $X$ be a normal variety. 
    A \emph{foliation} $\mathcal{F}$ on $X$ is a coherent subsheaf of the tangent sheaf $\mathcal{T}_X$ such that
\begin{itemize}
    \item $\mathcal{F}$ is saturated, that is, $\mathcal{T}_X/\mathcal{F}$ is torsion free and
    \item $\mathcal{F}$ is closed under Lie bracket.
\end{itemize}
\end{defn}
\begin{defn}
    Suppose $\mathcal{F}$ is a rank $r$ foliation on a normal variety $X$. 
    Notice that there exists an open embedding $j\colon X_0\hookrightarrow X$ such that $X_0$ is smooth, $\on{codim}(X\setminus X_0)\geq2$ and $\mathcal{F}$ is locally free on $X_0$, so $\wedge^r\mathcal{F}$ is an invertible sheaf on $X_0$. 
    We define the \emph{canonical divisor} of $\mathcal{F}$ to be any divisor $K_{\mathcal{F}}$ on $X$ such that $\mathcal{O}_X(-K_{\mathcal{F}})\cong j_*(\wedge^r\mathcal{F})$.
\end{defn}
\begin{defn}
    Let $X$ be a normal variety and $\mathcal{F}$ a rank $r$ foliation on $X$. 
    A subvariety $S \subseteq X$ is called $\mathcal{F}$-\emph{invariant} if for any open subset $U \subseteq X$ and any section $\partial \in \on{H}^{0}(U,\mathcal{F})$ we have:
    \[\partial(I_{S \cap U}) \subseteq I_{S\cap U}\]
    where $I_{S\cap U}$ denotes the ideal sheaf of $S \cap U$ in $U$. 
    If $\Delta \subseteq X$ is a prime divisor, one defines $\epsilon(\Delta)=1$ if $\Delta$ is invariant in the above sense and  $\epsilon(\Delta)=0$ otherwise.
\end{defn}
Let $f\colon X\dashrightarrow Y$ be a dominant map between normal varieties, $U\subseteq X_{\textnormal{smooth}}$, $V\subseteq Y_{\textnormal{smooth}}$ open subsets, such that $f\vert_U\colon U\to V$ is a smooth morphism. Then $f$ induces a foliation $\mathcal{G}$ on $X$, such that $\mathcal{G}_U$ is given by the tangent vectors of fibres of $f\vert_U$. For the precise definition of $\mathcal{G}$, we refer to \cite[Section 3.2]{druel2021codimension}.
\begin{defn}
    Let $\mathcal{F}$ be a foliation on a normal variety $X$. We say that $\mathcal{F}$ is \emph{algebraically integrable} if it is induced by some dominant map $f\colon X\dashrightarrow Y$. In this case we also call a divisor $D$ \emph{vertical} if it is $\mathcal{F}$-invariant, and \emph{horizontal} if it is not $\mathcal{F}$-invariant.
\end{defn}

We would like to mention that for algebraically integrable foliations, the best resolution we have is weakly semistable (see \cite[Theorem 2.1]{abramovich2000weak}). 
Thus, we cannot guarantee the existence of F-dlt modification. Instead, we use property ($*$)-modification. For the definition of property ($*$) and its corresponding modification, we refer to \cite[Section 3.5]{ambro2021positivity}.

\begin{defn}
    Let $X$ be a normal variety, $B$ be an effective $\mathbb{Q}$-divisor on $X$, and $\mathcal{F}$ be a foliation on $X$. 
    Then we call $(\cF,B)$ a \emph{foliated pair} (or $(X,\mathcal{F},B)$ a \emph{foliated triple}) if $K_{\mathcal{F}}+B$ is $\mathbb{Q}$-Cariter. Let $f\colon Y\to X$ be any birational morphism from some normal variety $Y$, let $\mathcal{G}$ be the induced foliation of $\mathcal{F}$ on $Y$, we will have
    \[K_{\mathcal{G}}+B_Y=f^*(K_{\mathcal{F}}+B)+\sum a(E_i,\cF,B)E_i\]
    where $B_Y$ is the strict transform of $B$ on $Y$. 
    We say that the triple $(X,\mathcal{F},B)$ is \emph{log canonical} if for any such model $Y$ over $X$, we have $a(E_i,\cF,B)\geq -\epsilon(E_i)$. 
    The triple $(X,\cF,B)$ is called \emph{strictly log canonical} if there is an exceptional divisor $E$ over $X$ such that $a(E,\cF,B)=-\epsilon(E)=-1$. 
    
    Let $V\subseteq X$ be a subvariety. 
    We say that $V$ is an \emph{lc} (resp. \emph{strictly lc} resp. \emph{non-lc}) \emph{centre} of $(X,\mathcal{F},B)$ if there exists a divisor $E$ over $X$ such that\begin{itemize}
        \item $a(E,\mathcal{F},B)=-\epsilon(E)$ (resp. $a(E,\mathcal{F},B)=-\epsilon(E)=-1$ resp. $a(E,\mathcal{F},B)<-\epsilon(E)$) and
        \item $\mathrm{centre}_X(E)=V$.
    \end{itemize}

    We denote by $\on{Nlc}(\cF,B)$ the union of all the non-lc centres of $(\cF,B)$. 
\end{defn}

\begin{defn}
    Let $X$ be a normal variety and $(\mathcal{F},B)$ be a foliated pair on $X$. 
    We say that $(\mathcal{F},B)$ is a \emph{Fano foliated pair} (or $(X,\mathcal{F},B)$ a \emph{Fano foliated triple}) if the following conditions are satisfied:
    \begin{enumerate}
        \item $X$ is klt,
        \item $-(K_\mathcal{F}+B)$ is ample, and
        \item $(\mathcal{F},B)$ is log canonical.
    \end{enumerate}
    We call $(X,\mathcal{F},B)$ \emph{weak Fano} if it satisfies (1) and (3) above, and $-(K_\mathcal{F}+B)$ is big and nef. 
\end{defn}
We recall the definition of complements for classical pairs. 
\begin{defn}
    Let $(X,B)$ be a pair. Then an $n$-complement is of the form $K_X+B^+$ such that
    \begin{itemize}
        \item $(X,B^+)$ is log canonical, 
        \item $n(K_X+B^+)\sim0$, and
        \item $nB^+\geq n\lfloor B\rfloor+\lfloor(n+1)\{B\}\rfloor$.
    \end{itemize}
\end{defn}
For the existence of $n$-complements for Fano pairs with bounded $n$, we refer to \cite[Theorem 1.10]{Antipluri}. 
Similarly, we define $n$-complements for foliations.
\begin{defn}
    Let $(X,\mathcal{F},B)$ be a foliated triple. Then an \emph{$n$-complement} is of the form $K_\mathcal{F}+B^+$ such that
    \begin{itemize}
        \item $(X,\mathcal{F},B^+)$ is log canonical, 
        \item $n(K_{\mathcal{F}}+B^+)\sim0$, and
        \item $nB^+\geq n\lfloor B\rfloor+\lfloor(n+1)\{B\}\rfloor$.
    \end{itemize}
\end{defn}

\begin{ex}[Unbounded log canonical Fano foliations]
Let $\bF_n$ be the Hirzebruch surface and $E$ be the negative section. 
Let $\cG_n$ be the foliation induced by the canonical fibration on $\bF_n$ and $\pi\colon \bF_n \to S_n$ be the contraction of $E$. 
Let $\cF_n=\pi_*\cG_n$ be the pushforward of $\cG_n$. 
Then we have the following:
\begin{enumerate}
    \item $(S_n,\cF_n)$ is strictly log canonical, 
    \item $-K_{\cF_n}$ is ample, and 
    \item $S_n$ has only one singularity, which is $\epsilon$-lc if and only if $\epsilon\geq \frac{2}{n}$. 
\end{enumerate}

Similarly, we can consider ruled surfaces $\bP(\cE)$ over curves of arbitrary genus and $\pi\colon \bP(\cE)\to S$ be the contraction of the negative section. 
Let $\cG$ be the foliation induced by the fibration and $\cF$ be the pushforward of $\cG$ by $\pi$. 
Then $\cF$ is strictly log canonical and $S$ has only one singularity which is possibly non-rational.
\end{ex}

\section{Fano foliations on klt surfaces}
In this section, we construct $1$-complements for Fano foliations of rank one on klt surfaces. 
\begin{prop}\label{cplt_surface}
    Let $\cF$ be a Fano foliation of rank one on a klt surface $X$. 
    Then $\cF$ has a $1$-complement, that is, there is an effective divisor $B$ such that the foliated pair $(\cF,B)$ is log canonical and $K_\cF+B\sim 0$. 
\end{prop}
\begin{proof}
    We divide the proof into the following steps, which will first state the claim, followed by the justifications:
    
    \begin{enumerate}[label=(Step \arabic*)]
        \item $\cF$ is algebraically integrable and strictly log canonical. 

        As $\cF$ is Fano, by bend and break, $\cF$ is algebraically integrable. 
        By \cite[Proposition 5.3]{araujo2013fano}, $\cF$ is strictly log canonical. 

        \item There is a birational morphism $\pi\colon W\to X$ extracting only non-invariant exceptional divisors $E_i$ over lc centre. 
        Moreover, the singularities of $W$ on $E_i$ are cyclic quotient singularities. 

        By \cite[Lemma 3.2]{araujo2013fano}, there is a fibration $\widetilde{f}\colon \widetilde{X}\to C$ parametrizing the general leaves of $\cF$. 
        Let $\widetilde{\pi}\colon \widetilde{X}\to X$ be the canonical morphism. 
        Note that $C\cong\bP^1$ as, over a strictly lc centre, there is a rational $\pi$-exceptional curve dominating $C$. 
        Now we run a $K_{\widetilde{X}}$-MMP over $\bP^1$ to contract invariant $\pi$-exceptional curves. 
        So it is also a $K_{\widetilde{X}}$-MMP over $X$. 
        Let $W$ be the output of this MMP with $\pi\colon W\to X$ and $f\colon W\to C$. 
        Note that, by \cite{chen2023log}, there is exactly one non-invariant divisor $E_i$ over each strictly lc centre $p_i$ and the singularities of $W$ on $E_i$ are all cyclic quotient singularities. 

        \item There is exactly one strictly lc centre $p_0$ on $X$ with $E_0$ the non-invariant exceptional divisor on $W$ over $p_0$. 
        Thus, the relative Picard number $\rho(W/X)=1$ and $\cF_{W}=\pi^{-1}\cF$ is induced by a morphism $f\colon W\to \bP^1$. 

        Let $F$ be a general fibre of $f$. 
        Then 
        \[0>K_\cF\cdot \pi_*F = \pi^*K_\cF\cdot F = (K_{\cF_W}+\sum_iE_i)\cdot F = -2 + \sum_i 1.\] 
        Thus, there is exactly one $p_i$, say $p_0$. 
        Therefore, $p_0\notin\on{Supp}\Delta$ by \cite[Lemma 2.2]{chen2022acc}
        
        \item For any fibre $F$ of $f$, if there is only one singularity of $W$ on $F$, then it is not a cyclic quotient singularity. 

        Suppose there is a fibre $F$ of $f$ on which there is only one singularities of $W$, say $p$, which is a cyclic quotient singularity. 
        Then let $E_1$, $\ldots$, $E_\ell$ be all prime exceptional divisors on the minimal resolution $g\colon Y\to X$ of $X$ over $p$ with $E_i\cdot E_{i+1}=1$ for $i=1$, $\ldots$, $\ell-1$. 
        Note that the intersection matrix $(E_i\cdot E_j)$ is negative definite. 
        
        Let $\widetilde{F}=g^*F+\sum_ia_iE_i$ with $\widetilde{F}\cdot E_1=1$. 
        Note that $0\leq\widetilde{F}\cdot E_j \leq (-\sum_iE_i)\cdot E_j$ for all $j$ and both inequalities are strict when $j=1$. 
        Thus, by \cite[Corollary 4.2]{kollar_mori_1998}, we have $a_i\in (-1,0)$ for all $i$.  
        Thus, $a_1 = F^2+a_1 = \widetilde{F}^2$ is an integer, which is impossible.

        \item For any fibre $F$ of $f$, if $W$ is smooth at $F\cap E_0$, then $W$ is smooth at all points on $F$. 

        Let $p\in F\setminus E_0$ be a point at which $W$ is singular. 
        Then $\cF$ is canonical at $p$ as there is exactly one strictly lc centre on $X$. 
        Note that \[0>K_\cF\cdot \pi_*F = \pi^*K_\cF\cdot F = (K_{\cF_W}+E_0)\cdot F = K_{\cF_W}\cdot F + 1\] 
        and thus $K_{\cF_W}\cdot F<-1$. 
        If $\cF$ is not terminal at $p$, then $K_{\cF_W}\cdot F\geq 1-2= -1$, which is impossible. 
        Thus, all singular points of $W$ on $F\setminus E_0$ are terminal singularities for $\cF_W$. 
        As terminal foliation singularities are supported at cyclic quotient singularities, by the previous step, there are at least two distinct singular points of $W$ on $F\setminus E_0$, then $K_{\cF_W}\cdot F\geq 2-2-\frac{1}{r_1}-\frac{1}{r_2}\geq -1$ where $r_1$, $r_2$ are two positive integer at least $2$, which is impossible. 

        \item There are at most three singular fibres of $f$, on which there are exactly two singular points of $W$ and one of them is on $E_0$. 
        Moreover, all singularities are cyclic quotient singularities. 

        By \cite{chen2023log}, as $X$ is klt at $p_0$, there are at most three singularities on $E_0$, say $q_i$, all of which are cyclic quotient singularities. 
        For each $q_i$, the fibre $F_i$ of $f$ through $q_i$ must have another singularities as $F_i^2=0$. 
        
        Moreover, by similar arguments above, there is exactly one other singularities, say $q'_i$, which is terminal foliation singularity. 

        \item When $W$ is smooth, there is a $1$-complement for $\cF$. 
        
        If $W$ is smooth, then it is a Hirzebruch surface $\bF_n$. 
        Let $C$ be a section of $f$ with $C\cap E_0=\varnothing$. 
        Thus, $(\cF_W,E_0+C)$ is a $1$-complement for $(\cF_W,E_0)$. 
        Pushing forward to $X$, we get a $1$-complement for $\cF$. 

        \item There is a $1$-complement for $\cF$ if $W$ is not smooth. 
        
        Suppose $W$ is not smooth, there are birational morphisms $\mu\colon Y\to W$ and $\nu\colon Y\to \bF_n$ where $\bF_n$ is a Hirzebruch surface. 
        Moreover, $\mu$ and $\nu$ are a sequence of at most three weighted blowups. 
        \[\xymatrix{ & Y \ar[ld]_-\mu \ar[rd]^-\nu & \\ W \ar[rd]_-{f} \ar@{-->}[rr]^-{h} & & \bF_n \ar[ld]^-{g} \\ & \bP^1. &}\]

        More precisely, let $C_0$ be the negative section for $\bF_n \to \bP^1$ and $L_i$ be a fibre for $\bF_n\to \bP^1$ such that $h_*^{-1}L_i$ is a singular fibre for $f$. 
        By Step 6, there are at most three $L_i$. 
        At each $L_i\cap C_0$, $\nu$ is a weighted blowup. 
        And $\mu$ contracts the strict transform of $L_i$ on $Y$. 

        Let $\cG$ be the foliation on $\bF_n$ induced by $g$ and $h_*(K_{\cF_W}+E_0) = K_{\cG}+C_0$ where $C_0$ is the negative section for $g$ and $h=\nu\circ\mu^{-1}$. 
        We have seen that there is a $1$-complement for $(\cG,C_0)$, say $(\cG,C_0+C)$ where $C$ is a section for $g$ with $C\cap C_0=\varnothing$. 
        Thus, pulling back to $W$ and then pushing forward to $X$, we get a $1$-complement $(\cF,\pi_*h_*^{-1}C)$ for $\cF$. 
    \end{enumerate}
\end{proof}

\section{Algebraically integrable Fano foliations}
In this section, we show the existence of $\bQ$-complements for Fano foliated triples $(X,\cF,B)$ with $\cF$ being algebraically integrable.\par
First, we prove adjunction and inversion of adjunction for invariant divisors in this case:
\begin{prop}[{$=$ \cite[Theorem 2.4.1]{chen2023minimal}, Adjunction}]\label{restriction}
    Let $(X,\mathcal{F},B)$ be a foliated triple. 
    Suppose $\mathcal{F}$ is algebraically integrable and $S\subseteq X$ a prime invariant divisor with normalisation $\nu\colon\Tilde{S}\to S$. 
    Assume that $(X,\mathcal{F},B)$ is log canonical, then $(\tilde{S},\mathcal{F}_{\tilde{S}},\on{Diff}_{\tilde{S}}(\mathcal{F},B))$ is log canonical.
\end{prop}
\begin{proof}
    We give a direct proof here, which is essentially the same as the one given in \cite{chen2023minimal}. 
    
    Let $f\colon (Y,\mathcal{G},E)\to (X,\mathcal{F},B)$ be a property $(*)$ modification of $(X,\mathcal{F},B)$. 
    Then we have
    \[K_{\mathcal{G}}+E=f^*(K_{\mathcal{F}}+B).\]
    Let $T$ be the strict transform of $S$ with normalisation $\tilde{T}\to T$. By taking the restriction to $\tilde{T}$, we have
    \[K_{\mathcal{G}_{\Tilde{T}}}+\on{Diff}_{\Tilde{T}}(\mathcal{G},E) = f\vert_{\Tilde{T}}^*(K_{\mathcal{F}_{\tilde{S}}}+\on{Diff}_{\tilde{S}}(\mathcal{F},B)).\]
    Since the left hand side is log canonical by \cite[Proposition 3.2]{ambro2021positivity}, hence the right hand side.
\end{proof}
\begin{lemma}\label{modi}
    Let $(X,\mathcal{F},B)$ be a Fano foliated triple and a divisor $\Delta\sim_\bQ-(K_{\mathcal{F}}+B)$. 
    Suppose $L\subseteq X$ is a non-lc centre of $(X,\mathcal{F},B+\Delta)$. 
    Then there exists a Property $(*)$ modification $\pi\colon(Y,\mathcal{G},\Tilde{B}+\tilde{\Delta}+E)\to (X,\mathcal{F},B+\Delta)$ such that 
    there exists a codimension one non-log canonical centre $W$ of $(Y,\mathcal{G},\tilde{\Delta}+\Tilde{B}+E)$ such that $\pi(W)=L$ where $E = \sum\epsilon(E_i)E_i$.
\end{lemma}
\begin{proof}
    We consider the local log canonical threshold
    \[t=\on{lct}_L(\mathcal{F},B;\Delta):=\sup\{\lambda\in\mathbb{R}\vert L\nsubseteq\on{Nlc}(\mathcal{F},B+\lambda\Delta)\}\in[0,1).\]
    Then $L$ is an lc centre of $(\mathcal{F},B+t\Delta)$. 
    By \cite[Theorem 3.10]{ambro2021positivity}, there exists a Property $(*)$ modification $\pi\colon(Y,\mathcal{G},\Tilde{B}+t\tilde{\Delta}+\sum\epsilon(E_i)E_i)\to (X,\mathcal{F},B+t\Delta)$ with a codimension one log canonical centre $W$ of $(Y,\mathcal{G},\Tilde{B}+t\tilde{\Delta}+\sum\epsilon(E_i)E_i))$ such that $\pi(W)=L$. 
    To be precise, we have
    \[K_{\mathcal{G}}+\Tilde{B}+t\tilde{\Delta}+\sum\epsilon(E_i)E_i+F^\prime=\pi^*(K_\mathcal{F}+B+t\Delta).\]
    Therefore we conclude
    \[K_{\mathcal{G}}+\Tilde{B}+\tilde{\Delta}+\sum\epsilon(E_i)E_i+(F^\prime+(1-t)(\pi^*(\Delta)-\tilde{\Delta}))=\pi^*(K_\mathcal{F}+B+\Delta).\]
    We denote the non-lc part by
    \[F:=F^\prime+(1-t)(\pi^*(\Delta)-\tilde{\Delta}).\]
    Since $(\mathcal{F},B)$ is log canonical but $(\mathcal{F},B+\Delta)$ is not log canonical along $L$, we have $L\subseteq\on{Supp}(\Delta)$. 
    Furthermore, as $\pi(W)=L$ we get
    \[\on{mult}_W(F)\geq\epsilon(W)+(1-t)\on{mult}_W(\pi^*(\Delta)-\tilde{\Delta})>\epsilon(W),\]
    hence the result.
\end{proof}
Next, we prove the following inversion of adjunction theorem for foliated pairs. We will need the following lemma:
\begin{lemma}\label{Sing}
    Let $\pi\colon Y\to Z$ be an equi-dimensional fibration between normal varieties with $\mathcal{F}$ the induced foliation on $Y$ such that $\dim Y=n$ and $\mathcal{F}$ is Gorenstein and has rank $r$. 
    Let $H_Z\subseteq Z$ be a prime divisor and assume that
    \begin{itemize}
        \item $S$ and $T$ are two irreducible components of $\pi^{-1}(H_Z)$, 
        \item $G\subseteq S\cap T$ is a component of $S\cap T$ and $\mathrm{codim}(G;S)=1$ such that $S$ and $T$ are both Cartier around $\eta_G$, and
        \item $\pi(G)=H_Z$.
    \end{itemize}
    Then $G\subseteq\on{Sing}(\mathcal{F})$.
\end{lemma}
\begin{proof}
    Since the statement is local around the generic point $\eta_G$, we can assume that:
    \begin{itemize}
        \item Both $Y$ and $Z$ are affine. 
        \item $Z$ and $H_Z$ are both smooth. In particular, we may assume that $\Omega_{Z}$ is locally free and generated by $\textnormal{d}z_1$, $\textnormal{d}z_2$, $\dots$, $\textnormal{d}z_{n-r}$ where $H_Z:=\{z_1=0\}$. 
        \item $\pi^\sharp(z_1)=x^{n_1}y^{n_2}f$ where $x$, $y$ are the local generators of $S$ and $T$ respectively, $n_1$, $n_2\geq1$ and $f\neq0$ a holomorphic function.
    \end{itemize}
    
    Now we consider the first fundamental exact sequence \cite[Theorem 25.1]{matsumura1970commutative}:
    \[\pi^*\Omega_{Z}\to\Omega_{Y}\to\Omega_{Y/Z}\to0.\]
    Notice that
    $\Omega_{Y/Z}\cong\Omega_{Y}/\langle\textnormal{d}(\pi^\sharp z_1),\textnormal{d}(\pi^\sharp z_2),\dots,\textnormal{d}(\pi^\sharp z_{n-r})\rangle$ 
    and that
    \[\textnormal{d}(\pi^\sharp z_1)=x^{n_1}\textnormal{d}(y^{n_2}f)+y^{n_2}\textnormal{d}(x^{n_1}f)\in\mathfrak{m}_{\eta_G}\Omega_{Y/\on{Spec}\bC}.\]
    Thus, we have $\dim_{k(\eta_G)}\Omega_{Y/Z}\otimes k(\eta_G)\geq r+1$ and hence, $\dim_{k(\eta_G)}\Omega_{Y/Z}^r\otimes k(\eta_G)>1$. 

    Now we consider the induced morphism $\varphi_r\colon\Omega_{Y/Z}^r\to\Omega_{Y/Z}^{[r]} = \mathcal{O}_Y(K_{\mathcal{F}})$. 
    Since $\mathcal{O}_Y(K_{\mathcal{F}})$ is locally free by assumption, $\varphi_r$ factors by some inclusion $\psi_r\colon\Omega_{Y/Z}^r/(\Omega_{Y/Z}^r)^{\textnormal{tor}}\hookrightarrow\mathcal{O}_Y(K_{\mathcal{F}})$: 
    \begin{center}
        \begin{tikzcd}
            \Omega_{Y/Z}^r\arrow[r,twoheadrightarrow]\arrow[dr,"\varphi_r",swap]&\Omega_{Y/Z}^r/(\Omega_{Y/Z}^r)^{\textnormal{tor}}\arrow[d,hookrightarrow,"\psi_r"]\\
            &\mathcal{O}_Y(K_{\mathcal{F}}).
        \end{tikzcd}
    \end{center}
    
    Since $\dim_{k(\eta_G)}(\Omega_{Y/Z}^r/(\Omega_{Y/Z}^r)^{\textnormal{tor}})\otimes k(\eta_G)=\dim_{k(\eta_G)}\Omega_{Y/Z}^r\otimes k(\eta_G)>1$, we have $\varphi_r$ is not surjective at $\eta_G$. 
    Then $\on{Im}(\Omega_{Y}^r\to\mathcal{O}_Y(K_{\mathcal{F}}))=\on{Im}(\Omega_{Y/Z}^r\xrightarrow{\varphi_r}\mathcal{O}_Y(K_{\mathcal{F}}))=\on{Im}\psi_r$ has cosupport containing $G$. 
    Thus, $\on{Im}(\Omega_{Y}^r(-K_\cF)\to\mathcal{O}_Y) \subseteq \cI_G$ and therefore, $G\subseteq\on{Sing}(\cF)$. 
\end{proof}

\begin{thm}[Inversion of adjunction]\label{inv}
    Let $(X,\mathcal{F},B)$ be a foliated triple with $\mathcal{F}$ being algebraically integrable. 
    Let $S\subseteq X$ be a prime invariant divisor with normalisation $\nu\colon S^\nu\to S$ such that $(S^\nu,\mathcal{F}_{S^\nu},\on{Diff}_{S^\nu}(\mathcal{F},B))$ is log canonical. Then
    \[\on{Nlc}(X,\mathcal{F},B)\cap S=\varnothing.\]
\end{thm}
\begin{proof}
    Suppose $\on{Nlc}(X,\mathcal{F},B)\cap S\neq\varnothing$ for the sake of contradiction. 

    By \cite[Theorem 3.10]{ambro2021positivity}, there is a Property $(*)$ modification $f\colon Y \to X$ such that $(Y,\cG,B_Y)$ satisfies Property $(*)$ with 
    \begin{equation}\label{crepant}
        K_{\mathcal{G}}+B_Y+F=f^*(K_\mathcal{F}+B)
    \end{equation}
    where $B_Y = f_*^{-1}B+\sum\epsilon(E_i)E_i$ and $F$ is effective and $f$-exceptional. 
    Let $\pi\colon Y\to Z$ be the morphism inducing $\mathcal{G}$ on $Y$. 
    Let $T:=f_*^{-1}(S)$ be the strict transform of $S$. 
    By assumption that $\on{Nlc}(X,\mathcal{F},B)\cap S\neq\varnothing$, we have $\on{Supp}(F)\cap T\neq\varnothing$. 
    Let $\rho\colon T^\rho\to T$ be the normalisation of $T$. 
    We have the following commutative diagram:
    \begin{center}
        \begin{tikzcd}
            T^\rho\arrow[d]\arrow[r,"\rho"]&T\arrow[r,hook]\arrow[d]&(Y,\mathcal{G},\tilde{B}+\tilde{E})\arrow[d,"f"]\arrow[r,"\pi"]&Z\\
            S^\nu\arrow[r,"\nu"]&S\arrow[r,hook]&(X,\mathcal{F},B)\arrow[ru,dashrightarrow]&
        \end{tikzcd}
    \end{center}
    Now we consider the adjunction of equation \ref{crepant} to $S$ and $T$:
    \begin{equation}\label{subcrepant}
        K_{\mathcal{G}_{T^\rho}}+\on{Diff}_{T^\rho}(\mathcal{G},B_Y)+F\vert_{T^\rho}=f^*(K_{\mathcal{F}_{S^\nu}}+\on{Diff}_{S^\nu}(\mathcal{F},B))
    \end{equation}
    Since we have the hypothesis that the restriction $(S^\nu,\mathcal{F}_{S^\nu},\on{Diff}(\mathcal{F},B))$ is log canonical, so is 
    \[(T^\rho,\mathcal{G}_{T^\rho},\on{Diff}_{T^\rho}(\mathcal{G},B_Y)+F\vert_{T^\rho}).\]
    In particular, we have therefore
    \[\on{mult}_G(\on{Diff}_{T^\rho}(\mathcal{G},B_Y)+F\vert_{T^\rho})\leq\epsilon(G)\]
    for any prime divisor $G\subseteq T^\rho$. 
    Since $\on{Supp}(F)\cap T\neq\varnothing$ and $Y$ is $\mathbb{Q}$-factorial, we may take $G$ to be a component of $\on{Supp}(F)\cap T$ such that $\on{codim}_T(G)=1$ and $G=F_1\cap T$ where $F_1\subseteq\on{Supp}(F)$. 
    Let $a=\on{mult}_{F_1}F$. 
    Note that $a>0$, $\on{mult}_G(F_1\vert_T)=a$, and 
    \begin{align*}
        \epsilon(G) 
        &\geq \on{mult}_G(\on{Diff}_{T^\rho}(\mathcal{G},B_Y)+F\vert_{T^\rho}) \\
        &\geq \on{mult}_G\on{Diff}_{T^\rho}(\mathcal{G}) + \on{mult}_G((B_Y+F)\vert_T) \\
        &\geq \on{mult}_G\on{Diff}_{T^\rho}(\mathcal{G}) + \on{mult}_{F_1}(B_Y+F) \\
        &\geq \on{mult}_G\on{Diff}_{T^\rho}(\mathcal{G}) + \epsilon(F_1)+a \\
        &\geq a > 0.
    \end{align*}
    Thus, $\epsilon(G)=1$, $\epsilon(F_1)=0$, and 
    \begin{equation}\label{mult_G_less_1}
        \on{mult}_G\on{Diff}_{T^\rho}(\cG)\leq 1-a<1.
    \end{equation}
    Therefore, $G\subseteq\on{Sing}(\mathcal{G})$ by Lemma \ref{Sing}. 
    
    Let $p\colon(Y^\prime,\mathcal{G}^\prime)\to(Y,\mathcal{G})$ be the quasi-\'{e}tale index one cover of $K_{\mathcal{G}}$, $T$ and, $F$ around $\eta_G$, with $T^\prime$ the normalisiotn of $p_*^{-1}T$, $G^\prime:=p_*^{-1}G$. 
    Let $e$ be the ramification index of $p$ along $G^\prime$. 
    By \cite[Lemma 3.4]{druel2021codimension} (see also \cite[Proposition 2.2]{cascini2021mmp}), we have
    \[\on{mult}_{G^\prime}\on{Diff}_{T^\prime}(\mathcal{G}^\prime) = e\on{mult}_{G}(\on{Diff}_{T}\mathcal{G})-(e-1).\]
    As $G\subseteq\on{Sing}(\mathcal{G})$ and $\cG^\prime$ is Cartier around $\eta_{G^\prime}$, we have $\on{mult}_{G^\prime}\on{Diff}_{T^\prime}(\mathcal{G}^\prime)\geq1$ and thus $\on{mult}_{G}\on{Diff}_{T}(\mathcal{G})\geq 1$, which contradicts the inequality \ref{mult_G_less_1}.  
\end{proof}
Next, we prove that the log canonical on a general fibre implies global log canonical property:
\begin{prop}\label{horizofNlc}
    Let $(Y,\mathcal{G},\tilde{B})$ be a foliated triple where $\mathcal{G}$ is induced by a dominant rational map to a normal variety $Z_Y$ and such that $K_{\mathcal{G}}+\tilde{B} \sim_{\mathbb{Q}}0$. Then $\on{Nlc}(Y,\mathcal{G},\tilde{B})$ dominates $Z_Y$.
\end{prop}
\begin{proof}
    Let $\pi\colon (X,\mathcal{F},B) \to (Y,\mathcal{G},\tilde{B})$ be a Property $(*)$ modification of $(Y,\cG,\tilde{B}_{\textnormal{h}})$ where $\tilde{B}_{\textnormal{h}}$ is the horizontal part of $\tilde{B}$, $B=\pi_*^{-1}\tilde{B}_{\textnormal{h}}$, and $\cF$ is a foliation induced by a morphism $X\to Z$. 
    We write 
    \[K_\cF+B+G=\pi^*(K_\cG+\tilde{B})\]
    where $G$ is an effective divisor. 
    Note that $\on{Nlc}(X,\mathcal{F},B+G)$ is supported in $G$ and by construction, $K_{\mathcal{F}}+B+G = \pi^*(K_{\mathcal{G}}+\tilde{B})\sim_{\mathbb{Q}} 0$. 
    As $K_\cF+B\sim_{\bQ} -G$ is not pseudo-effective over $Z$, by \cite[Thereom 1.4]{LMX2024minimal}, we can run a $(K_{\mathcal{F}}+B)$-MMP$/Z$ with scaling of an ample$/Z$ $\bR$-divisor which terminates with $(W,\cF_W,B_W)$ a Mori fibre space$/Z$ of $(X,\cF,B)$ where $B_W$ and $\mathcal{F}_W$ denote the pushfoward of the divisor $B$ on $X$ and the foliation $\cF$, respectively.   
    Let $\phi\colon X\dashrightarrow W$ be the birational map induced by the MMP and $\psi\colon W\to T$ be the contraction introduced by the Mori fibre space structure. 

    Let $G_W$ be the pushfoward of the divisor $G$ on $X$. 
    We observe that, by construction, $G_W$ is ample over $T$ and hence, in particular horizontal over $T$. 
    Furthermore, since $K_{\mathcal{F}_{W}}+B_{W}+G_{W} = \phi_{*}(K_{\mathcal{F}}+B+G)\sim_{\mathbb{Q}}0$, $-(K_{\mathcal{F}_W}+B_W)\sim_{\mathbb{Q}}G_W$ is ample over $T$ as well. 
    We label some morphisms in the following commutative diagram: 
    \begin{center}
    \begin{tikzcd}
    X \arrow[r,dotted, "\phi"] \arrow[d]
     & W \arrow[d, "\psi"] \arrow[dl, "\nu"] \\
    Z 
     & T\arrow[l, "\eta"]
    \end{tikzcd}
    \end{center}
    As $\nu$ is surjective, so is $\eta$, which implies that $G_W$ dominates $Z$. 
    Therefore, $G$ also dominates $Z$. 
    The claim follows upon observing that $\pi_{*}(\on{Nlc}(X,\mathcal{F},B+G)) \subseteq \on{Nlc}(Y,\mathcal{G},\tilde{B})$.
\end{proof}
\begin{rmk}
Notice that this proposition together with Lemma \ref{inv} proves that lifts of complements from vertical divisors $S$ to the ambient space remain log canonical. 
\end{rmk}

Now we are able to show the existence of $\mathbb{Q}$-complements for algebraically integrable Fano foliations.

\begin{thm}\label{Qcom}
    Let $(X,\mathcal{F},B)$ be an lc Fano foliated triple with $\cF$ algebraically integrable. 
    Then there exists a $\bQ$-complement. 
    That is, there exists an effective divisor $\Delta\sim_\bQ-(K_{\mathcal{F}}+B)$ such that $(X,\mathcal{F},B+\Delta)$ is log canonical.
\end{thm}
\begin{proof}
    We will proceed by induction on the corank $c$ of $\cF$. 
    When $c=0$, the existence of $\bQ$-complement is given by \cite[Theorem 1.7]{Antipluri}. 

    Now we assume the theorem holds when the corank is $c-1$. 
    Let $f\colon (Y,\mathcal{G},B_Y+E)\to(X,\mathcal{F},B)$ be a Property $(*)$ modification of $(X,\mathcal{F},B)$ such that $\mathcal{G}$ is induced by $\pi\colon Y\to Z$ and 
    \[K_{\mathcal{G}}+B_Y+E=f^*(K_{\mathcal{F}}+B)\]
    where $B_Y=f_*^{-1}B$ is the strict transform of $B$ on $Y$ and $E=\sum\epsilon(E_i)E_i$ is the sum of all non-invariant $f$-exceptional divisors. 
    We fix a general prime divisor $D_Z\subseteq Z$ such that $F:=\pi^{-1}(D_Z)$ is normal and irreducible. 
    Let $G:=f(F)$ be the image. 
    Notice that $G$ is not necessarily normal. 
    Therefore, we denote the normalisation by $\nu\colon G^\nu\to G$ and the normality of $F$ naturally gives us an induced morphism $g^\nu\colon F\to G^\nu$.
    \begin{center}
        \begin{tikzcd}
            &F\arrow[d,"g"]\arrow[dl,"g^\nu",swap]\arrow[r,"\iota_F",hook]&(Y,\mathcal{G},B_Y+E)\arrow[r,"\pi"]\arrow[d,"f"]&Z\\
            G^\nu\arrow[r,"\nu"]&G\arrow[r,"\iota_G",hook]&(X,\mathcal{F},B)\arrow[ru,dashrightarrow]&
        \end{tikzcd}
    \end{center}
    
    By adjunction, $(\cF_{G^\nu}, \on{Diff}_{G^\nu}(\cF,B))$ is an lc Fano foliated pair of corank $c-1$. 
    Thus by induction hypothesis, there exists an effective divisor $\Delta_{G^\nu}$ on $G^\nu$ such that $K_{\cF_{G^\nu}}+\on{Diff}_{G^\nu}(\cF,B)+\Delta_{G^\nu}\sim_\bQ 0$ and $(\cF_{G^\nu}, \on{Diff}_{G^\nu}(\cF,B)+\Delta_{G^\nu})$ is lc. 
    Let $m_0\in\bN$ such that $m_0\Delta_{G^\nu}\sim -m_0(K_{\cF_{G^\nu}}+\on{Diff}_{G^\nu}(\cF,B))$ and $L := -m(K_\cF+B)$ where $m$ is sufficiently divisible by $m_0$. 
    Then 
    \[L\vert_{G^\nu} =(\iota_G\circ\nu)^*L = -m(K_{\cF_{G^\nu}}+\on{Diff}_{G^\nu}(\cF,B))\] 
    and thus, $m\Delta_{G^\nu}\in |L\vert_{G^\nu}|$. 
    Therefore, there is a non-zero section $\alpha_{G^\nu}\in\on{H}^0(G^\nu,\cO_{G^\nu}(L\vert_{G^\nu}))$ such that $m\Delta_{G^\nu} = \on{Div}(\alpha_{G^\nu})+L\vert_{G^\nu}$. 

    Now we consider the following short exact sequence induced by the normalisation $\nu$ and twisted by the line bundle $\mathcal{O}_G(L\vert_G)$:
    \[0 \to \cO_{G}(L\vert_G) \to \nu_*\cO_{G^\nu}(L\vert_G) \to (\nu_*\cO_{G^\nu}/\cO_G)(L\vert_G) \to 0\]
    and thus, we have the following exact sequence:
    \[0 \to \on{H}^0(G,\cO_G(L\vert_G)) \xrightarrow{j} \on{H}^0(G,\nu_*\cO_{G^\nu}(L\vert_G)) \xrightarrow{\beta} \on{H}^0(G,(\nu_*\cO_{G^\nu}/\cO_G)(L\vert_G)).\]

    Note that $\alpha_{G^\nu}\in\on{H}^0(G,\nu_*\cO_{G^\nu}(L\vert_G))$ and $E\vert_F$ is a union of some lc centres of $(\cG_F,\on{Diff}_F(\cG,B_Y+E))$ and thus, $g^\nu(E\vert_F)$ is a unioin of some lc centres of $(\cF_{G^\nu},\on{Diff}_{G^\nu}(\cF,B))$. 
    Thus, no component of $\Delta_{G^\nu}$ is contained in $g^\nu(E\vert_F)$. 
    Hence, the image $\beta(\alpha_{G^\nu})$ is $0$ and therefore, there is an $\alpha_G\in\on{H}^0(G,\cO_G(L\vert_G))$ such that $j(\alpha_G) = \alpha_{G^\nu}$. 

    We then consider the following exact sequence on $X$:
    \[0 \to \cO_X(L-G) \to \cO_X(L) \to \cO_G(L) \to 0.\]
    We recall that $L=-m(K_\cF+B)$ is an ample divisor. 
    Then by Serre vanishing theorem, we have $\on{H}^1(X,\cO_X(L-G))=0$ for $m$ sufficiently large and thus, the following morphism is surjective:
    \[\on{H}^0(X,\cO_X(L))\twoheadrightarrow \on{H}^0(G,\cO_{G}(L)).\]
    Let $\alpha\in \on{H}^0(X,\cO_X(L))$ be a lifting of $\alpha_G$. 
    Then we have an effective divisor $D :=\on{Div}(\alpha)+L$ and thus $K_\cF+B+\Delta\sim_\bQ 0$ where $\Delta = \frac{1}{m}D$. 

    Note that $K_\cG+B_Y+E+f^*\Delta = f^*(K_\cF+B+\Delta)$. 
    Taking adjunction to $F$, we have
    \begin{align*}
        K_{\cG_F}+\on{Diff}_{F}(\mathcal{G},B_Y+E+f^*\Delta) 
        &= (g^\nu)^*(K_{\cF_{G^\nu}}+\on{Diff}_{G^\nu}(\mathcal{F},B+\Delta)) \\
        &= (g^\nu)^*(K_{\cF_{G^\nu}}+\on{Diff}_{G^\nu}(\mathcal{F},B)+\Delta_{G^\nu}).
    \end{align*}
    Since $(\cF_{G^\nu},\on{Diff}_{G^\nu}(\cF,B)+\Delta_{G^\nu})$ is log canonical, so is $(\cG_F,\on{Diff}_{F}(\cG,B_Y+E+f^*\Delta))$. 
    
    According to the inversion of adjunction (Theorem \ref{inv}), we have 
    \[\on{Nlc}(\cG,B_Y+E+f^*\Delta)\cap F=\varnothing\]
    and thus $\on{Nlc}(\cG,B_Y+E+f^*\Delta)$ cannot be horizontal as $F:=\pi^{-1}(D_Z)$ is irreducible. 
    Hence, by Proposition \ref{horizofNlc}, $\on{Nlc}(\cG,B_Y+E+f^*\Delta)=\varnothing$. 
    Therefore, $(Y,\cG,B_Y+E+f^*\Delta)$ is log canonical and so is $(X,\cF,B+\Delta)$.
\end{proof}

We now make note of a couple of immediate corollaries of the proof of the above proposition, which are not made further use of in the present paper, however appear to be of independent interest. 
The first can be seen as log-extension of a classification result obtained in \cite{araujo2014fano}, the subsequent bounds on the volumes of  the foliated anti-canonical divisor may be of interest in light of a possible foliated version of the BAB conjecture for algebarically integrable Fano foliations some of the authors are working towards.

We notice that Proposition \ref{horizofNlc} provides us with an alternative algebraic proof of \cite[Proposition 3.14]{araujo2014fano}:
\begin{prop}
    Let $(X,\mathcal{F},B)$ be an lc Fano foliated triple, then $\mathcal{F}$ is not induced by a morphism.
\end{prop}
\begin{proof}
    Suppose, for the sake of contradiction, that $\mathcal{F}$ is induced by a morphism $\pi\colon X\to Z$. 
    We will then proceed the proof by induction on $c$. 
    When $c=1$, we have $\dim Z=1$. 
    We fix a general point $z_0\in Z$ and $F_0:=\pi^{-1}(z_0)$. 

    Let $z_1\in Z\setminus\{z_0\}$ be another general point and $F_1:=\pi^{-1}(z_1)$ such that $F_1$ is irreducible and normal. 
    Note that $(F_1,\on{Diff}_{F_1}(\cF,B))$ is a Fano pair, there exists a $\bQ$-complement $D_{z_1}$ such that $K_{F_1}+\on{Diff}_{F_1}(\cF,B)+D_{z_1}\sim_\bQ 0$. 
    Since $-(K_{\mathcal{F}}+B)$ is ample, there exists an $m\gg0$ such that $L:=-m(K_{\mathcal{F}}+B)-F_0$ is very ample. 
    
    Note that $mD_{z_1}\sim L\vert_{F_1}$. 
    We consider the following short exact sequence on $X$:
    \[0\to\cO_X(L-F_1)\to\mathcal{O}_X(L)\to\mathcal{O}_{F_1}(L)\to0.\]
    By Serre vanishing theorem, we have the induced surjective morphism:
    \[\on{H}^0(X,\mathcal{O}_X(L))\twoheadrightarrow \on{H}^0(F_z,\mathcal{O}_{F_z}(L)).\]
    Let $mD\in \on{H}^0(X,\mathcal{O}_X(L))$ be a lifting of $mD_z$. 
    Then the triple $(X,\mathcal{F},B+D+\frac{1}{m}\Delta)$ has the following properties:
    \begin{itemize}
        \item $K_{\mathcal{F}}+B+D+\frac{1}{m}F_0\sim_\mathbb{Q}0$.
        \item $(F_1,\on{Diff}_{F_1}(\mathcal{F},B+D+\frac{1}{m}F_0))=(F_1,\on{Diff}_{F_1}(\mathcal{F},B)+D_{z_1})$ is log canonical since $z_0\neq z_1$. 
        \item $(X,\mathcal{F},B+D+\frac{1}{m}F_0)$ is obviously not log canonical since $F_0>0$ is a vertical divisor.
    \end{itemize}
    Thus, we get a contradiction with Proposition \ref{horizofNlc}. 

    Now we suppose the theorem holds for the cases when corank $\leq c-1$. 
    We fix a general Cartier divisor $H_Z>0$ on $Z$ and $S:=\pi^*H_Z$. 
    Then $(S,\cF_S,\on{Diff}_S(\cF,B))$ is an lc Fano triple with corank $c-1$ and $\cF_S$ is a foliation induced by the morphism $\pi\vert_S\colon S\to H_Z$, which contradicts the induction hypothesis. 
\end{proof}

Last, as an addendum, we prove a version of connectedness principle for foliations induced by morphisms:
\begin{lemma}\label{conlem}
    Let $(\mathcal{F},\Delta)$ be a foliated pair on a normal variety $X$ such that
    \begin{enumerate}
        \item $-(K_{\mathcal{F}}+\Delta)$ is big and nef;
        \item $(\mathcal{F},\Delta)$ is log canonical;
        \item $\mathcal{F}$ is induced by a morphism $g\colon X\to Z$.
    \end{enumerate}
    Then $\lfloor\Delta\rfloor$ is connected. 
\end{lemma}
\begin{proof}
    We prove the result by induction on dimension of $X$. Suppose the result is true for dimension $n-1$. 
    
    Let $H$ be a general vertical divisor with respect to $g$, that is, $H=g^*A$ for some general ample divisor $A$. 
    Let $m\colon X^\prime\to X$ be a log resolution of $(\mathcal{F},\Delta)$ such that the strict transform $H^\prime:=m_*^{-1}H$ is also smooth. 
    Let $\mathcal{F}^\prime$ be the induced foliation on $X^\prime$. 
    We consider the adjunction with respect to $H^\prime$:
    \[(K_{\mathcal{F}}+\Delta)\vert_{H^\prime}=K_{\mathcal{F}_{H^\prime}}+\Theta=K_{\mathcal{F}_{H^\prime}}+\Delta\vert_{H^\prime}+\Theta_0.\]
    According to \cite[Proposition 3.2]{ambro2021positivity}, we have $\on{Supp}(\lfloor\Theta_0\rfloor)\subseteq\on{Sing}(\mathcal{F}^\prime)\cap H^\prime$. 
    However since $H=g^*A$ for some general $A$ and $\on{codim}(\on{Sing}(\mathcal{F}))\geq2$, by moving $H$ in an infinitesimal neighbourhood in $\mathbb{P}(\on{H}^0(Z,\mathcal{O}_Z(A)))$ if necessary, we can assume that we have $\on{Supp}(\lfloor\Theta_0\rfloor)=0$. 
    Here we use \cite[Proposition 3.2]{ambro2021positivity} again to see that $(\mathcal{F}_H,\Theta)$ is a weak Fano lc foliation pair. 
    Hence we can use our induction assumption to get that $\lfloor\Theta\rfloor=\lfloor\Delta\vert_H\rfloor$ is connected. 
    Therefore, $\lfloor\Delta\rfloor$ is connected. 
    
    Now we only have to show the base case. 
    However, it is already proven in \cite[Theorem 17.4]{kollar1992flips}.
\end{proof}

\bibliographystyle{plain} 
\bibliography{Fano}

\end{document}